\documentclass{article}

\usepackage{amsmath}
\usepackage{url}

\usepackage{amssymb}
\usepackage{amsmath}
\usepackage{amsthm}
\usepackage{amsfonts}
\usepackage{mathrsfs}
\usepackage{makeidx}
\usepackage{multicol}

\begin{document}

\title{A characterization of rational functions}
\author{Bao Qin Li}
\date{}
\maketitle

\begin{abstract}
We give an elementary characterization of rational functions among
meromorphic functions in the complex plane.
\end{abstract}



\vskip.5in

\noindent While entire functions (i.e., complex-valued functions
differentiable in the complex plane ${\bf C}$) generalize
polynomials, meromorphic functions in ${\bf C}$ (i.e.,
complex-valued functions differentiable in ${\bf C}$ except at
poles) generalize rational functions. These are far-reaching
generalizations and have been important objects in the function
theory of complex variables and various applications. Characterizing
polynomials among entire functions and rational functions among
meromorphic functions is a natural goal. Different characterizations
in different contexts can be found in the literature (see \cite{Fu}, p.52, \cite{GO}, p.34, \cite{Hu}, p.143, etc.)


Given a polynomial $f=a_mz^m+\cdots+a_1z+a_0$ of degree $m,$ it is
immediate to see that
$\lim\limits_{z\rightarrow\infty}z\frac{f'(z)}{ f(z)}=m.$ This limit
condition turns out to be also sufficient for an entire function to
be a polynomial of degree $m$ (see \cite{Li}). It is natural to ask
if this condition can be pushed over to characterizing rational
functions among meromorphic functions. This short note gives an
affirmative answer to this question. The characterization, which is elementary and of a neat form, does not seem to have
been observed before and we thought it might be of some use to present it to the reader. 

Given a rational function $f$, we can write it as $f={P\over Q}$
with $P$ and $Q$ co-prime polynomials (i.e., $P$ and $Q$ have no
common polynomial factors). The number $\deg P-\deg Q$, the
difference between the degrees of $P$ and $Q$, is, by abuse of
language, called the divisor of $f$. When $f$ reduces to a
polynomial, the divisor is just the degree of the polynomial. We
then have the following

\bigskip
\noindent {\bf Theorem 1 (Rational function characterization).} {\it
A nonzero meromorphic function $f$ is a rational function of divisor
$d$ if and only if $\lim\limits_{z\rightarrow\infty}z\frac{f'(z)}{
f(z)}=d$. } \noindent
\begin{proof} The necessary condition is quite obvious. Write $f={P\over Q}$, where $P$ and $Q$ are co-prime polynomials of degree $m$ and $n$,
respectively. By definition, $d=m-n$. It is easy to check that
${f'\over f}={P'\over P}-{Q'\over Q}$. Thus,
\begin{eqnarray*}
& &\lim_{z\rightarrow\infty}z\frac{f'(z)}{
f(z)}=\lim_{z\rightarrow\infty}\{z\frac{P'(z)}{P(z)}-z\frac{Q'(z)}{Q(z)}\}=m-n=d.
\end{eqnarray*}

To prove sufficiency, we first show that $f$ must have finitely many
zeros and poles. In fact, by the given limit condition,
$|z\frac{f'(z)}{ f(z)}|\le M$ for some $M>0$ and sufficiently large
$|z|$, say for $|z|\ge R>0$ . Note that a zero or a pole $z_0\not=0$
of $f$ must be a pole of $z\frac{f'(z)}{ f(z)}$, which then tends to
infinity as $z\to z_0$. Thus, $z_0$ must lie inside the disc
$|z|<R$. This shows that all the zeros and poles lie inside this
disc, and thus $f$ can only have finitely many zeros and poles (Otherwise, 
these zeros/poles would have a finite limit point, which implies that $f$ is identically zero or has a non-pole singularity, a contradiction).
Assume that $f$ has $m$ zeros and $n$ poles (counting multiplicities/orders). We can write
\begin{eqnarray}\label{function}
f=\frac{p(z)}{q(z)}e^{g(z)},
\end{eqnarray}
where $g$ is an entire function, $p(z)=(z-a_1)\cdots (z-a_m), q(z)=(z-b_1)\cdots (z-b_n)$, and $a_j$'s (resp. $b_j$'s) are zeros (resp. poles) of $f$ (repeated according to their multiplicites/orders). Take $r>R$, then we have that
$m-n={1\over 2\pi i}\int_{|z|=r}{f'(z)\over f(z)}dz,$ where the
circle is traversed counterclockwise. This follows from the Argument
Principle or from Cauchy's theorem directly using (\ref{function}):
\begin{eqnarray*}
& &{1\over 2\pi i}\int_{|z|=r}{f'(z)\over f(z)}dz\\
& &=\sum_{k=1}^m{1\over 2\pi i}\int_{|z|=r}{dz\over
z-a_k}-\sum_{k=1}^m{1\over 2\pi i}\int_{|z|=r}{dz\over
z-b_k}\\
& &+{1\over 2\pi i}\int_{|z|=r}g'(z)dz=m-n.
\end{eqnarray*}
But,
\begin{eqnarray*}
& &{1\over 2\pi i}\int_{|z|=r}\frac{f'(z)}{f(z)}dz={1\over 2\pi
}\int_0^{2\pi}re^{i\theta}\frac{f'(re^{i\theta})}{f(re^{i\theta})}d\theta\rightarrow
d
\end{eqnarray*}
as $|z|=r\rightarrow \infty$, since
$$re^{i\theta}\frac{f'(re^{i\theta})}{f(re^{i\theta})}=z\frac{f'(z)}{f(z)} \rightarrow d$$
as $|z|=r\rightarrow \infty$ uniformly in $\theta$ by the given
condition of the theorem. Thus, we obtain that
\begin{eqnarray}\label{divisor}
m-n=d.
\end{eqnarray}
It is easy to check that
$$\lim_{z\rightarrow\infty}(z\frac{p'(z)}{p(z)}-z\frac{q'(z)}{q(z)})=m-n$$ and
$$\frac{f'(z)}{f(z)}=\frac{p'(z)}{p(z)}-\frac{q'(z)}{q(z)}+g'(z),$$
which implies, by the given limit condition and (\ref{divisor}),
that
\begin{eqnarray}\label{zero}
& &\lim_{z\rightarrow\infty}zg'(z)\\
&
&=\lim_{z\rightarrow\infty}(z\frac{f'(z)}{f(z)}-\bigl(z\frac{p'(z)}{p(z)}-z\frac{q'(z)}{q(z)})\bigr)=d-(m-n)=d-d=0.
\end{eqnarray}
Thus, the entire function $zg'(z)$ is bounded and thus, by Liouville's
theorem, must be constant, which is then identically $0$ by (\ref{zero}). Hence, 
$g'(z)\equiv 0$, i.e, $g$ is constant, from which it follows that
$f=\frac{p(z)}{q(z)}e^{g(z)}$ is a rational function. This completes
the proof. \end{proof}

From the sufficiency of Theorem 1 (cf. the proof), we see that {\it if $\lim\limits_{z\rightarrow\infty}z\frac{f'(z)}{ f(z)}=d$ then $f$ is a rational function with $m$ zeros and $n$ poles
satisfying that $m-n=d$}, which generalizes the Fundamental Theorem
of Algebra: A nonzero polynomial $f$ of degree $d$ clearly satisfies
that $\lim\limits_{z\rightarrow\infty}z\frac{f'(z)}{ f(z)}=d$ and
thus has $d$ zeros (counting multiplicities).

\bigskip
\noindent
\textit{Department of Mathematics and Statistics\\
Florida International University\\
Miami, FL 33199 USA\\
libaoqin@fiu.edu}

\end{document}